\documentclass{amsart}
\begin{document}
\title{Mathematical Education}
\author{William P. Thurston}
\maketitle

\noindent {\bf Note:}
This article originally appeared in the \emph{Notices of the AMS} {\bf 37}
(1990), 844--850.  Copyright 1990 William P. Thurston.
\vskip 1pc

Mathematics education is in an unacceptable state. Despite much
popular  attention to this fact, real change is slow.  

Policymakers often do not comprehend the nature of mathematics or of 
mathematics education. The `reforms' being implemented in different school 
systems are often in opposite directions. This phenomenon is a sign that 
what we need is a better understanding of the problems, not just the recognition
that they exist and that they are important.

I am optimistic that our nation will find solutions to these problems. Problems
arising from failure of understanding are curable. We do not lack for 
dedication, resources, or intelligence: we lack direction.
\section{Symptoms}

There are many symptoms of the problems in mathematics education. The number 
of undergraduate mathematics majors is about half what it was 15 years ago. The
number of U.S. graduate students is less than half what it was 15 years ago, 
although foreign students have taken up some of the slack.  The performance of
our students at all levels, as measured by standardized tests, is below that 
of other industrialized countries.

The typical response of American adults, on meeting a mathematician, is one of
dismay.  They apologetically recall the last mathematics course they took, 
which is usually the one where they lost their grip on the subject matter.

Many companies recognize lack of mathematical competence as a major problem in
their workforce.  Technology has removed the need for elementary mathematical 
facility in a few jobs, such as ringing up hamburger orders, but it has also 
eliminated many unskilled jobs, such as assembly-line work, while creating 
many others requiring considerable mathematical sophistication.

Students in college mathematics courses are unresponsive. They are afraid to 
speculate and afraid to reach into themselves for ideas. When one visits 
classrooms at different grade levels, one sees a dramatic decline in 
liveliness and spontaneity with age.  One gets the impression that the natural
interest and curiosity of young children in mathematics has been weeded out.

In most places, it is harder to get students to do homework or to study outside
the classroom than it was twenty years ago.

Even those college students who have been successful in the high school 
mathematics curriculum, including calculus, have a narrow base of knowledge of
mathematics.


\section{ Stratification and Compartmentalization}

A major source of problems, and a major barrier to solutions of the problems, 
is the stratification and compartmentalization of the immense mathematical 
education enterprise, from kindergarten through graduate school.  In 
particular, there is very little communication between high school and college
teachers of mathematics.  There is also less than optimal communication between
the people who are involved with curricular reform at the college level and 
research mathematicians.  This splintering partly results from the existence 
of several different professional associations, including the American 
Mathematical Society (AMS), which primarily represents mathematical research;
the Mathematical Association of America (MAA), which primarily represents
undergraduate mathematics; and the National Council of Teachers of Mathematics
(NCTM), which primarily represents the most committed high school mathematics
teachers.  Additional organizations represent two-year colleges, applied 
mathematics, statistics, and operations research, as well as computer science.

The membership of the AMS and the MAA overlap by only about $1/3$. Most 
members of the AMS are not even aware of the existence of the NCTM. Conversely,
few members of the NCTM are aware of the AMS.  Very few academic mathematicians
have any contact with school mathematics teachers except through their
children. Most teachers are trained at institutions with few if any research 
mathematicians.  A number of universities are strong both in mathematical 
research and teacher training, but, even at those universities, most research 
mathematicians are not involved in teacher training.

Even more severe is the compartmentalization due to the division of primary 
and secondary education into many individual school districts and schools where
the real decisions are made, and the division of college education into many
different academic departments where the real decisions are made. Many people 
are thinking about the problems and making individual efforts to solve them.  
Many of these local efforts are quite successful.  Unfortunately, there is far 
too little coordination and sharing of the insights gained.

\section{ Mathematics is a Tall Subject}

One feature of mathematics which requires special care in education 
is its `height,' that is, the extent to which concepts build on previous 
concepts.  Reasoning in mathematics can be very clear and certain, and, once 
a principle is established, it can be relied on.  This means it is possible 
to build conceptual structures which are at once very tall, very reliable, 
and extremely powerful.  

The structure is not like a tree, but more like a scaffolding, with many 
interconnected supports.  Once the scaffolding is solidly in place, it is not 
hard to build it higher, but it is impossible to build a layer before previous
layers are in place.

Difficulties arise because students taking a particular course are in different
stages of mastery of the earlier learning.  They also tend to be secretive 
about exactly what they know and what they don't know. For instance, 
many calculus students don't correctly add fractions, at least not 
in symbolic form: the typical mistake is that 
$$ \frac ab + \frac cd = \frac{a + c}{b + d} $$
much simpler than
$$ \frac{a d + c b }{b d} $$
However, students feel guilty that they are shaky on addition of fractions and
are slow to admit it and ask how and why it works, even to themselves.

Addition of fractions is a very boring topic to someone who already knows it, 
but it is an essential skill for algebra, which in turn is essential for 
calculus. It is not so hard, when talking with students individually, to find 
out what parts of the structure need shoring up and to deal with those parts 
individually.  But it is quite difficult to find a level of teaching which is 
comprehensible and at the same time interesting to an entire class with 
heterogeneous background.


\section{ Mathematics is a Broad Subject}

Mathematics is also very broad.  There are many subjects that are 
never discussed in the mainline curriculum which culminates in calculus. 
The subjects that are discussed have many interesting side-branches 
that are never explored.

In my parents' generation (during the 1940s), the standard first college 
mathematics course was college algebra.  Soon afterward, the standard first 
college course was calculus, until the early 1960s, when calculus became 
standard for the best high school mathematics students.  By now first year 
calculus has largely migrated to high school in affluent school districts, so 
that most of the better mathematics and science students at our best 
universities have already taken calculus before they arrive. At Princeton, for
instance, two-thirds of entering students placed out of at least one semester
of calculus last year.

The acceleration of the curriculum has had its cost: there has been an 
accompanying trend to prune away side topics.  When I was in high school, for 
instance, it was standard to study solid geometry and spherical geometry along
with plane geometry.  These topics have long been abandoned. The shape of the 
mathematics education of a typical student is tall and spindly.  It reaches a 
certain height above which its base can support no more growth, and there it 
halts or fails.

These two trends (lengthening and narrowing of school mathematics) have been 
hastened by the growing reliance on standardized tests. Standardized tests 
are designed to cover topics on the most standard curriculum: if only half the 
students study some topic, it is unfair to ask about it on a standardized test.
This is not so bad as long as tests are used in a disinterested way as one of 
several {\it devices for measurement}.   Instead, higher test scores are often
treated as the {\it goal}. Legislators, newspapers and parents put 
superintendents and school boards under pressure, superintendents and school 
boards put principals under pressure, principals put teachers under pressure;
and teachers put students under pressure to raise scores.  The sad result is 
that many mathematics courses are specifically designed to raise scores on 
some standardized test.

We don't diagnose pneumonia with only a thermometer, and we don't attempt to
cure it by putting ice in a patient's mouth.  We should take a similarly
enlightened attitude toward testing in mathematics education.

The long-range objectives of mathematics education would be better served if 
the tall shape of mathematics were de-emphasized, by moving away from a 
standard sequence to a more diversified curriculum with more topics that start
closer to the ground.  There have been some trends in this direction, such as 
courses in finite mathematics and in probability, but there is room 
for much more.


\section{ Mathematics is Intuitive and Real}

Students commonly lose touch with the reality and the intuitive nature of
mathematics.  From kindergarten through high school, they often have teachers
who are uncomfortable with anything off the beaten path. Young children come up with 
many ingenious devices to work out mathematical questions, but teachers usually
discourage any nonconventional approach---partly because it is not easy to 
understand what a child is thinking or trying to say, and the teachers don't 
catch on, partly because the teachers think it's not okay to use an alternative
method or explanation.

By the time students are in college, they are inhibited from thinking for 
themselves and from admitting out loud what they are thinking.  Instead, they 
try to figure out what routines they are supposed to learn.  When there is any
departure in class from the syllabus or the text, someone invariably asks 
whether it's going to be on the test. 

Unless mathematics makes a real connection to people, they are unlikely
ever to think about it or use it once they have completed a course.


\section{ Precocity and Competition}

Along with the emphasis on tests has come an emphasis on precocity and 
acceleration in mathematics.   It is relatively easy for a bright student to 
work through the mathematics curriculum far more quickly than the usual pace.

There are several problems associated with precocity. People who skip ahead 
in the curriculum often have gaps in their background which only show up later.
At that point, the person may be too embarrassed to admit the gap and tries to
fake understanding.  This regularly leads to disastrous results.

Another problem is that precocious students get the idea that the reward is in 
being `ahead' of others in the same age group, rather than in the quality of 
learning and thinking.  With a lifetime to learn, this is a shortsighted 
attitude.  By the time they are $25$ or $30$, they are judged not by 
precociousness but on the quality of work.  It is often a big letdown to 
precocious students when others who are talented but not so precocious catch 
up, and they become one among many. The problem is compounded by parents in 
affluent school districts who often push their children to advance as quickly 
as possible through the curriculum, before they are really ready.

A third problem associated with precociousness is the social problem. Younger 
students are often well able to handle mathematics classes intellectually 
without being able to fit in socially with the group of students taking them.  

Related to precociousness is the popular tendency to think of mathematics as
a race or as an athletic competition.  There are widespread high school math 
leagues: teams from regional high schools meet periodically and are given 
several problems, with an hour or so to solve them.

There are also state, national and international competitions. These 
competitions are fun, interesting, and educationally effective for the people 
who are successful in them. But they also have a downside.  The competitions 
reinforce the notion that either you `have good math genes', or you do not.
They put an emphasis on being quick, at the expense of being deep and 
thoughtful.  They emphasize questions which are puzzles with some hidden 
trick, rather than more realistic problems where a systematic and persistent 
approach is important.

This discourages many people who are not as quick or as practiced, but might 
be good at working through problems when they have the time to think through 
them. Some of the best performers on the contests do become good 
mathematicians, but there are also many top mathematicians who were not so 
good on contest math.  Quickness is helpful in mathematics, but it is only 
one of the qualities which is helpful. For people who do not become 
mathematicians, the skills of contest math are probably even less relevant.

These contests are a bit like spelling bees.  There is some connection 
between good spelling and good writing, but the winner of the state spelling 
bee does not necessarily have the talent to become a good writer, and some 
fine writers are not good spellers. If there was a popular confusion between 
good spelling and good writing, many potential writers would be unnecessarily 
discouraged.

I think the answer to these problems is to build a system which exploits the 
breadth of mathematics, by allowing quicker students to work through the 
material in greater depth and to take excursions into related topics, before 
racing ahead of their age group.

\section{ Mystery and Mastery}

Mathematics is amazingly compressible: you may struggle a long time, step
by step, to work through some process or idea from several approaches.  But
once you really understand it and have the mental perspective to see it
as a whole, there is often a tremendous mental compression.  You can file it
away, recall it quickly and completely when you need it, and use it
as just one step in some other mental process.  The insight that goes
with this compression is one of the real joys of mathematics.

After mastering mathematical concepts, even after great effort,
it becomes very hard to put oneself
back in the frame of mind of someone to whom they are mysterious.

I remember as a child, in fifth grade, coming to the amazing (to me)
realization that the answer to $134$ divided by $29$ is $134/29$ (and so forth).
What a tremendous labor-saving device!  To me, `134 divided by 29'
meant a certain tedious chore, while  $134/29$ was an object
with no implicit work.  I went excitedly
to my father to explain my major discovery.   He told me that of course
this is so, $a/b$ and $a$ divided by $b$ are just synonyms.  
To him it was just a small variation in notation.

One of my students wrote about visiting an elementary school and being
asked to tutor a child in subtracting fractions.  He was
startled and sobered to see
how much is involved in learning this skill for the first time,
a skill which had condensed to a triviality in his mind.

Mathematics is full of this kind of thing, on all levels. It never stops.

The hard-earned and powerful tools which are available almost unconsciously to
mathematicians, but not to students, make it hard for mathematicians to learn
from their students. This puts a psychological barrier in the way of listening
fully to students.

It is important in teaching mathematics to work hard to overcome this barrier
and to get out of the way enough to give students the chance to work things out
for themselves.


\section{ Competence and Intimidation}

Similarly, students at more advanced levels know many things which less
advanced students don't yet know. It is very intimidating to hear others
casually toss around words and phrases as if any educated person should know
them, when you haven't the foggiest  idea what they're talking about.  Less
advanced students have trouble realizing that they will (or would) also learn
these theories and their associated vocabulary readily when the time comes and
afterwards use them casually and naturally.  I  remember many occasions when I
felt intimidated by mathematical words and concepts before I understood them: 
negative, decimal, long division, infinity, algebra, variable, equation,
calculus, integration, differentiation, manifold, vector, tensor, sheaf,
spectrum, {\it etc.}  It took me a long time before I caught on to the pattern
and developed some immunity.

Teachers also are frequently intimidated about mathematics. High school
teachers are often timid about approaching college and university teachers.  
They also question, with some justice, whether university professors are in
touch with the problems they have to deal with.  There is so little general
contact between those who teach mathematics in high schools, colleges, and
universities that few professors know much about the educational problems in
high school or elementary school. Elementary school teachers are often quite
insecure about their grasp of mathematics and timid about approaching anyone.


\section{ The Problems and Solutions are the Same}

There is much in common about the problems---and the solutions---of
mathematics education at different levels, from kindergarten through
graduate school.  The failure to communicate is a real loss, and the
potential gain from opening of two-way communications is great.

Over the past two years, I have met many people involved in mathematical
education at all levels, partly as a member of the MSEB (Mathematical Sciences
Education Board), a national board for mathematics education, composed of
people from diverse backgrounds in mathematics and education, including a
number of teachers. I have learned a lot.

During the spring semester of 1990, John Conway, Peter Doyle, and I
organized a new course (`Geometry and the Imagination') at Princeton
which borrowed a good deal from the ideas I learned.  We taught as
a team, shunning lectures and emphasizing group discussions among
students.  We emphasized manipulables, cooperative learning, and
problem solving.  We asked students to keep journals and to write out
their ideas in good and complete English. Each student did a major
project for the course.  These ideas are all borrowed
from ideas current in K-12 mathematics education, focused in the NCTM
curriculum standards.

To culminate, we held a `geometry fair': like a science fair but with a popcorn
machine and without prizes.  It was great fun.

The course came alive, qualitatively more than any course we had
taught before. Students learned a lot of mathematics and solved problems
we wouldn't have dared ask in a conventional college class.

One topic we discussed was mathematical education.    Students were given an
earlier draft of this essay and wrote 70 thoughtful essays based on their own
experiences. This essay has benefitted considerably from their comments.


\section{ Socially Excluded Groups}

Why do so few women become mathematicians, and so few of the non-Asian
minorities?

I am convinced that the poverty of the school mathematics teaching and
curriculum has a lot to do with it.   The way mathematics is taught in school
does not address the real goals of a mathematical education. It is very hard to
get a sense of the depth, liveliness, power, and breadth of mathematics from
any ordinary experience with mathematics in school.  I believe that most
students who really master the subject matter, and eventually become
scientists, mathematicians, computer programmers, {\it etc.},  are those who
have some other channel for learning mathematics, outside the classroom. 
Sometimes it is the home, sometimes it is books, sometimes it is an unusual
teacher, but often it is the `nerd' social subgroup in school. When I was in
high school, I certainly belonged to this subgroup (although the name `nerd'
was not yet current), and I appreciated it very much.  However, it is a very
different matter for white and Asian males to join this social subgroup than
for women, Hispanics or blacks.

Since channels of learning outside the classroom are currently dominant in
mathematics education for the top group of students, improvement of the general
quality of mathematics within the classroom therefore should act as an
equalizer, particularly helping blacks, Hispanics, and women.  For people who
are not already tuned in to the mathematical style of discourse, it is
especially important to teach in a way that is not watered down, but that
begins from a person's real experience.

Intimidation also has a lot to do with it.  The emphasis on precocity, high
test scores, and competition work to amplify the small differences which arise
from other sources.


\section{ Goals and Standardized Tests}

What is mathematics education good for?

\subsection{ Mathematics in life}

First, mathematics is a basic tool of everyday life.  When coffee comes in $13$
oz. packages and $16$ oz. cans, can you take that information in stride
(walking slowly by the shelf), along with the prices and the prominent signs
claiming `contains 23\% more' on the cans, to help decide which you'd prefer to
buy? In buying a new car with various gimmicks in financing, rebates, and
features, do you understand what is going on?  If most people did, the gimmicks
would be pointless.

Second, mathematical reasoning is an important part of informed citizenship. Can
you understand the reasoning behind studies of health risks from various
substances, and can you judge how important they are? When listening to
politicians, can you and do you use your reckoning powers to help decide how
important some statistic is, and what it means? Can you measure and calculate
adequately for simple sewing and carpentry?  Can you plan a budget? When you
see graphs in newspapers and magazines, do you understand what they mean, and
are you aware of the several devices frequently used either to dramatize or to
play down a certain trend?

Third, mathematics is a tool needed for many jobs in the infrastructure of our
increasingly complex and technological society.  These uses are pervasive and
varied.  The dental technician, the fax repairperson, the fast food manager,
the real estate agent, the computer consultant, the bookkeeper and the banker,
the nurse and the lawyer, all need a certain proficiency with mathematics in
their jobs.

Fourth, mathematics is intensively used (and sometimes abused) in most branches
of science. Much of theoretical science really {\it is} mathematics. 
Statistics is one of the most common uses of mathematics.  Many scientists use
the widespread computerized statistical packages, which alleviate the need for
computation. However, people who use statistical packages are often shaky in
their understanding of the basic principles involved and often apply
statistical tests or graphical displays inappropriately.


\subsection{ Mathematics is alive}

To me, these utilitarian goals are important, but secondary.  Mathematics has a
remarkable beauty, power, and coherence, more than we could have ever expected.
It is always changing, as we turn new corners and discover new delights and
unexpected connections with old familiar grounds.  The changes are rapid,
because of the solidity of the kind of reasoning involved in mathematics.

Mathematics is like a flight of fancy, but one in which the fanciful turns out
to be real and to have been present all along. Doing mathematics has the feel
of fanciful invention, but it is really a process of sharpening our perception
so that we discover patterns that are everywhere around. In his famous apology
for mathematics, G.H. Hardy praised number theory for its purity, its
abstraction, and the self-evident impossibility of ever putting it to practical
use.  Now this very subject is applied very widely, particularly for encoding
and decoding communications.

My experience as a mathematician has convinced me that the aesthetic goals and
the utilitarian goals for mathematics turn out, in the end, to be quite close. 
Our aesthetic instincts draw us to mathematics of a certain depth and
connectivity.  The very depth and beauty of the patterns makes them likely to
be manifested, in unexpected ways, in other parts of mathematics, science, and
the world.

To share in the delight and the intellectual experience of mathematics---to
fly where before we walked---that is the goal of a mathematical education.


\subsection{ Testing and ``accountability''}

Unfortunately the goals of school mathematics have become incredibly narrow,
much narrower even than the first set of utilitarian goals listed above, let
alone the others.  It is popular lately for politicians and the public to
demand ``accountability'' from the school systems.  This would be great, except
that educational accounting is usually based on narrowly-focused
multiple-choice tests.

It is as if students were considered to have mastered Shakespeare if they could
pass a timed vocabulary test in Elizabethan English, or that they had learned
to write when they could correctly choose the grammatical form of a sentence
from four possibilities.

The state and regional boards of education these days hand out a laundry list
of skills which students are supposed to know at a certain age, rather than a
curriculum: horizontal addition versus vertical addition, addition of 2 digit
numbers to 2 digit numbers with a 2 digit answer versus addition of 2 digit
numbers to 2 digit numbers with a 3 digit answer, {\it etc.}   

A front-page article in the New York Times Metropolitan section of July 24,
1990 contrasted the elementary schools in two similar difficult neighborhoods
of Brooklyn. The first was a `successful' school, with two reading lessons a
day, the second lesson in `test-taking skills' and practice for the
standardized reading test. In this school, 80.5\% scored at or above grade
level.  The other school was an `unsuccessful' school where they prepare for
the test for `only' 3 months. In that school, 36.4\% score at or above grade
level.

The reporter cited an example of how the principal sets the tone in the
`successful' school:



\begin{quote}
She is not satisfied with just the right answers;  she wants the right steps
along the way.  In one fourth-grade class, she noticed that pigtailed Keanda
Snagg had made a wild, though accurate, stab at a problem requiring her to
average which of two stores had lower prices.

\vskip5pt

``It looks like you were going to do it without doing the work,'' she told
Keanda.

\vskip5pt

As she watched Keanda go through the calculations, she stressed fundamentals
beyond arithmetic, like putting numbers in clearly ordered columns.

\end{quote}

\vskip7pt

I can't tell which school is actually more successful without seeing them for
myself, but one thing I know:  neither the test scores  nor the cited incident
are demonstrations of greater success.

\section{ Thinking and Rote}

Narrow goals are stultifying.

People are much smarter when they can use their full intellect and when they
can relate what they are learning to situations or phenomena which are real to
them.

The natural reaction, when someone is having trouble understanding what you are
explaining, is to break up the explanation into smaller pieces and explain the
pieces one by one.  This tends not to work, so you back up even further and
fill in even more details.

But human minds do not work like computers: it is harder, not easier, to
understand something broken down into all the precise little rules  than to
grasp it as a whole.  It is very hard for a person to read a computer assembly
language program and figure out what it is about.  A computer reads and
executes it in the blink of an eye.  But the most powerful computer in the
world is not clever enough to drive a car safely, or control a stroll along the
sidewalk, or come up with an interesting mathematical discovery.

Studying mathematics one rule at a time is like studying a language by first
memorizing the vocabulary and the detailed linguistic rules, then building 
phrases and sentences, and only afterwards learning to read, write, and
converse.  Native speakers of a language are not aware of the linguistic rules:
they assimilate the language by focusing  on a higher level, and absorbing the
rules and patterns subconsciously. The rules and patterns are much harder for
people to learn explicitly than is the language itself.  In fact, the
tremendous and so far unsuccessful attempts to teach languages to computers
demonstrate that nobody can yet describe a language  adequately by precise
rules.

It is better not to teach a topic at all than to attempt teaching it in tiny
rules and bits.

\section{ Answers and Questions}

People appreciate and catch on to a mathematical theory much better after they
have first grappled for themselves with the questions the theory is designed to
answer.

There is a natural tendency, in teaching mathematics, to use the logical order
and to explain all the techniques and answers before bringing up the examples
and the questions, on the supposition that the students will be equipped with
all the techniques necessary to answer them when they arise. 

It is better to keep interesting unanswered questions and unexplained examples
in the air, whether or not the students, the teachers, or anybody is yet ready
to answer them.  The best psychological order for a subject in mathematics is
often quite different from the most efficient logical order.

As mathematicians, we know that there will never be an end to unanswered
questions. In contrast, students generally perceive mathematics as something
which is already cut and dried---they have just not gotten  very far in
digesting it.

We should present mathematics to our students in a way which is at once more
interesting and more like the real situations where students will encounter it
in their lives---with no guaranteed answer.


\section{ What Can We Do?}

This depends on who we are.

In our compartmentalized system, it is hard for a single organization or
individual to do very much to affect the overall system directly.  But
addressing the local situation will indirectly influence the global situation.
I will address the question from the point of view of college and university
mathematicians.

First, college and university mathematics departments should develop courses
which can give students a fresh chance in mathematics. Remedial courses are
widespread, but their success is limited: going over the same material one more
time is tedious and boring, whether you understood it the first time or you
didn't.  There is a built in handicap to enthusiasm and spontaneity.

Instead, there should be more courses available to freshmen and nonmajors which
exploit some of the breadth of mathematics, to permit starting near the ground
level without a lot of repetition of topics that students have already heard. 
For instance, elementary courses in topology, number theory, symmetry and group
theory, probability, finite mathematics, algebraic geometry, dynamical systems
(chaos), computer graphics and linear algebra, projective geometry and
perspective drawing, hyperbolic geometry, and mathematical logic can meet this
criterion.

Second, we should work to create better channels of communication between the
compartments of the educational system. We need to find devices so that the
educational accomplishments of professors are visible within the profession,
not merely within the classroom or within the department.  We need to find
vehicles for exchange of interesting ideas between different departments: for
instance, exchange visits between directors of undergraduate studies and chairs
of departments.  We should visit  each other's classrooms.  We need more talks
and special sessions related to education at our professional meetings, and
more prizes for educational accomplishments.

The newsletter, {\it Undergraduate Mathematics Education  Trends}, and the MAA
newsletter, {\it Focus}, are two such vehicles, and there are several other
publications which carry articles on undergraduate education. The annual
chairman's coloquium organized by the Board on Mathematical Sciences of the
National Research Council is another vehicle of communication, although it has
a broader agenda than education. Still, these existing channels are very small
compared with what we could establish.

Even more important and more difficult is the creation of channels for
communication between the strata: most important for colleges and universities
is communication between high school, college and university mathematics
departments.  
This communication must be two-way: college and university
professors can learn a lot about how to teach from school teachers. The MSEB is
one such channel, along with the state mathematics coalitions that they have
helped to stimulate, but what we need is a much more massive exchange.

How can the senior professors, who are at the top of a system which is clearly
not doing such a great job, presume to teach their juniors how to do better? 
The graduate students and the junior faculty often do a better job at teaching
mathematics than the senior faculty, who have sometimes become resigned to the
dismal situation, settled into a routine, and given up on trying any new
initiatives.  Even when they do a pretty good job in their own classrooms,
against the odds, they do not usually get involved in improving the overall
system.  

Often other professors are suspicious of the professor who does take an
interest in education. They tend to assume that research is the only activity
which really matters and that turning to education is a sign of failure in
research.   Senior professors sometimes explicitly advise junior faculty not to
waste too much energy on teaching, or they will never be promoted.

We must recognize that there are many different ways that we can make important
contributions to society and to our institutions. It is dumb to measure
mathematicians against the single scale of research.  Education is an important
and challenging endeavor, which many people engage in by choice, not necessity. 
We should judge them by what they accomplish, not by what they might have
accomplished if they spent their time and energy elsewhere.

What urgently needs to change is the system of professional rewards. We need
something better than the current situation within university mathematics
departments where there is lip service to the importance of teaching, but, when
it comes to the crunch of hiring and tenure decisions, teaching and service
count only in the marginal cases where the candidates cannot be differentiated
by the quality of research.

People are socially motivated.  As we discuss education with each other, we put
more energy into it, and it becomes more important to us.  The academic culture
{\it can} change, and it has changed.  The process of change is mostly an
informal one (what you talk about at lunch), not controlled by organizational
decisions. But when the time is ripe, as I believe it is now for mathematical
education, a little nudging by organizations can help stimulate a huge change.

The needed reforms will take place through collegial, cooperative efforts. Good
mathematical ideas spread very rapidly through informal channels in the
mathematical community.  As we turn more of our attention to education, good
educational ideas will also spread rapidly.

\end{document}